\documentclass[11pt,a4paper]{amsart}
\usepackage{amsfonts,amsmath,amssymb,amsthm,mathtools}
\usepackage[noadjust]{cite}

\usepackage{todonotes}

\numberwithin{equation}{section}

\DeclareMathOperator{\Ker}{Ker}
\DeclareMathOperator{\dist}{dist}
\DeclareMathOperator{\Dom}{Dom}
\DeclareMathOperator{\spec}{spec}
\DeclareMathOperator{\diag}{diag}

\DeclarePairedDelimiter{\abs}{|}{|}
\DeclarePairedDelimiter{\norm}{\lVert}{\rVert}

\DeclareSymbolFont{SY}{U}{psy}{m}{n}
\DeclareMathSymbol{\emptyset}{\mathord}{SY}{'306}

\newcommand{\CC}{\mathbb{C}}
\newcommand{\NN}{\mathbb{N}}

\newcommand{\RR}{\mathbb{R}}

\newcommand{\cH}{{\mathcal H}}

\newtheorem{theorem}{Theorem}[section]{\bf}{\it}
\newtheorem{lemma}[theorem]{Lemma}{\bf}{\it}
\newtheorem{corollary}[theorem]{Corollary}{\bf}{\it}
\newtheorem{proposition}[theorem]{Proposition}{\bf}{\it}
\newtheorem{remark}[theorem]{Remark}{\it}{\rm}
\newtheorem{example}[theorem]{Example}{\it}{\rm}

\title[]{Protecting points from operator pencils}
\date{}

\author[A.~Seelmann]{Albrecht Seelmann}
\address{A.~Seelmann,
Fakult\"at f\"ur Mathematik, Technische Univer\-si\-t\"at Dortmund,
D-44221 Dortmund, Germany}
\email{albrecht.seelmann@mathematik.tu-dortmund.de}

\author[M. T\"aufer]{Matthias T\"aufer}
\address{M. T\"aufer,
School of Mathematical Sciences, Queen Mary University of London, London, E1 4NS, United Kingdom
}
\email{m.taeufer@qmul.ac.uk}

\author[K.\ Veseli\'c]{Kre\v simir Veseli\'c}
\address{K. Veseli\'c,
Fakult\" at f\" ur Mathematik und Informatik, Fernuniversit\" at Hagen, D-58084 Hagen, Germany}
\email{kresimir.veselic@fernuni-hagen.de}

\begin{document}

\begin{abstract}
We classify all sets of the form $\bigcup_{t\in\mathbb{R}}\mathrm{spec}(A+tB)$ where $A$ and $B$ are self-adjoint operators and $B$
is bounded, non-negative, and non-zero. We show that these sets are exactly the complements of discrete subsets of $\mathbb{R}$,
that is, of at most countable subsets of $\mathbb{R}$ that contain none of their accumulation points.
\end{abstract}

\maketitle

\section{Introduction and main results}

We study the union of spectra $\spec( A + t B )$, $t \in \RR$, where $A$ and $B$ are self-adjoint operators in a (complex) Hilbert
space, and $B$ is bounded and non-negative. It is rather easy to cook up examples where this union is not the whole real line,
that is, there can exist \emph{protected points}, see, e.g., Example~\ref{example_2x2} below. The objective of this note is to
classify all possible sets of such protected points. Our main result reads as follows.

\begin{theorem}
\label{thm}
 Let $A$ and $B$ be self-adjoint operators on the same Hilbert space, and suppose that $B$ is bounded, non-negative, and non-zero.
 Then, the set $\RR \setminus \bigcup_{t\in\RR}\spec(A+tB)$ is discrete, that is, it is at most countable and
 contains none of its accumulation points. In particular, the union of these spectra is dense.

 Conversely, for every discrete $P \subset \RR$, there exist self-adjoint operators $A$ and $B$ with $B$ non-negative and
 bounded such that $P = \RR \setminus \bigcup_{t\in\RR}\spec(A+tB)$.
\end{theorem}

Theorem~\ref{thm} can be applied to obtain information on the unperturbed operator $A$ if it is known that $ \spec (A + tB) $ does
not vary with $t$.

\begin{corollary}[{See also~\cite{Ves18}}]
Let $A$ and $B$ be self-adjoint operators on the same Hilbert space, and suppose that $B$ is bounded and non-negative. If the
spectrum of $A + tB$ is independent of $t$, that is, $\spec(A + tB) = \spec(A)$ for all $t \in \RR$, then $B = 0$ or
$\spec(A) = \RR$.
\label{cor}
\end{corollary}

In fact, Corollary~\ref{cor} was the starting point of our investigation and finds an application in the context of the
Klein-Gordon equation in~\cite{Ves18} featuring the matrix operator
\begin{equation*}\label{eq:K-G}
    \mathbf{K} 
    =
    \begin{pmatrix}
        0      & \sqrt{p^2 +1}  \\
        \sqrt{p^2 +1}   &    2x \\
    \end{pmatrix},
    \quad
    p = -i \frac{\mathrm{d}}{\mathrm{d}x},
\end{equation*}
which has the property
\begin{equation} 
    \label{eq:homogeneous}
    e^{-ipt}\mathbf{K}e^{ipt} = \mathbf{K} -
    2t
    \begin{pmatrix}
        0 & 0  \\
        0  &  I \\
    \end{pmatrix}.
\end{equation}
Identity~\eqref{eq:homogeneous} resembles the definition of a \emph{homogeneous operator}. 
These are operators $A$ such that for every $t \in \RR$ the operator $A + t I$ is unitarily equivalent to $A$.
It is known that homogeneous operators have absolutely continuous spectrum on the whole real line~\cite{Tomsic-74} and a quantum
mechanical instance of such a phenomenon has been produced in (7.19) of~\cite{EnssV-83}. It would be interesting to study under
which assumptions it is possible to weaken the notion of homegeneity in~\cite{Tomsic-74} to the situation of the corollary above
and prove absolute continuity of the spectrum of $A$, and thus of $\mathbf{K}$ in~\eqref{eq:homogeneous}.

Also recall that the~\emph{spectrum of a linear operator pencil} is defined as
\begin{equation*}
 \spec(A,B)
 =
 \left\{ 
  \mu \in \CC 
  \colon
  0 \in \spec(A - \mu B)
 \right\},
\end{equation*}
see for instance~\cite{Marcus-88}. In this context, our considerations imply that for all sets $P$ as in Theorem~\ref{thm} there
exist self-adjoint operators $A$ and $B$ such that $\spec(A - \lambda, B)$ is empty for all $\lambda \in P$,
cf.~Proposition~\ref{prop} below.

In~\cite{GesztesyNN-15}, a related situation of linear operator pencils is studied. However, rather than \emph{protected points},
the authors investigate the set of parameters $t$ for which a given $\lambda \in \RR$ is in the \emph{point spectrum} of
$A + t B$. They review a result from~\cite{NabokoNS-13} that states that this set has Lebesgue measure zero and show that
sign-definiteness of $B$ is crucial for this to hold.
One easily sees that sign-definiteness of $B$ also cannot be dropped in Theorem~\ref{thm}. 
In fact, for indefinite $B$, complete intervals can be protected, as can be seen from the following example of $2 \times 2$
matrices 
(cf.~also Example~\ref{example_2x2} below):
\begin{equation*}
 A = \begin{pmatrix} 1 & 0\\ 0 & -1 \end{pmatrix},
 \quad
 B = \begin{pmatrix} 0 & 1\\ 1 & 0 \end{pmatrix}.
\end{equation*}
With these choices we have $\spec( A+tB ) = \{\pm\sqrt{1+t^2}\}$, so that the complement of the union
$\bigcup_{t\in\RR} \spec(A+tB)$ is equal to $(-1,1)$. This example is an instance of a more general result that states that a gap
in the spectrum of a self-adjoint operator is preserved under an off-diagonal perturbation, see Theorem~2.1 in~\cite{AdamjanL-95}
and also Theorem~8.1 in~\cite{DavisK-70}. Accordingly, also Corollary~\ref{cor} is no longer valid for general perturbations $B$,
that is, one can find self-adjoint operators $A$ and $B \neq 0$ with $B$ bounded and $\spec(A+tB) = \spec(A) \subsetneq \RR$ for
all $t\in\RR$.

\section{Proofs}

The first assertion of Theorem~\ref{thm} will follow from Lemmas~\ref{lem:count} and~\ref{lem:accumulation_points} below, the
second one from Lemma~\ref{lem:cyclic}. The proof of the latter is constructive. Corollary~\ref{cor} will then follow from the
first part of the theorem (cf.~also Corollary~\ref{cor2} below) and the fact the spectrum of $A$ is always closed.
The core of our considerations is the following result.

\begin{proposition}\label{prop}
 Let $A$ and $B$ be as in Corollary~\ref{cor}. Then, the following
 are equivalent:
 \begin{enumerate}
  \renewcommand{\theenumi}{\roman{enumi}}
  \item $0$ belongs to each resolvent set $\rho(A+tB)$, $t\in\RR$;
  \item $0\in\rho(A)$ and $BA^{-1}B=0$.
 \end{enumerate}
 In this case, it holds that $0\in\rho(A+zB)$ for all $z \in \CC$ with
 \begin{equation}
  \label{eq:inverse}
  (A+zB)^{-1} = A^{-1} - zA^{-1}BA^{-1}.
 \end{equation}
 
 \begin{proof}
  Suppose that $0 \in \rho(A)$. We first observe that
  \begin{equation}\label{eq:spec}
   \spec(BA^{-1})\setminus\{0\} = \spec(B^{1/2}A^{-1}B^{1/2}) \setminus \{0\} \subset \RR,
  \end{equation}
  where for the last inclusion we have used that $B^{1/2}A^{-1}B^{1/2}$ is self-adjoint. Moreover, we have for all
  $z \in \CC \setminus \{0\}$ that
  \begin{equation}\label{eq:res}
   A+zB = \bigl(I+zBA^{-1}\bigr)A = z\Bigl(\frac{1}{z}I+BA^{-1}\Bigr)A.
  \end{equation}
  Hence, $0\in\rho(A+zB)$ holds if and only if $0\in\rho(I/z+BA^{-1})$. Identity~\eqref{eq:spec} shows that the latter holds indeed
  for all $z \in \CC \setminus \RR$.

  (i)$\Rightarrow$(ii). By hypothesis, we have $0 \in \rho(A+zB)$ for all $z \in \RR$, whence, according to the considerations
  above, $0 \in \rho(I/z+BA^{-1})$ for all $z \in \CC \setminus \{0\}$. This yields $\spec(BA^{-1}) = \{0\}$, so that
  $\spec(B^{1/2}A^{-1}B^{1/2})=\{0\}$ by~\eqref{eq:spec}. From the self-adjointness of $B^{1/2}A^{-1}B^{1/2}$ we obtain
  $B^{1/2}A^{-1}B^{1/2}=0$, and, in particular, $BA^{-1}B=0$.

  (ii)$\Rightarrow$(i). We have $(BA^{-1})^2=0$. For each $z \in \CC$, the inverse of the operator $I+zBA^{-1}$ is therefore given
  by $I-zBA^{-1}$, and from~\eqref{eq:res} we conclude that
  \begin{equation*}
   (A+zB)^{-1} = A^{-1}(I-zBA^{-1})=A^{-1} - zA^{-1}BA^{-1}.
   \qedhere
  \end{equation*}
 \end{proof}%
\end{proposition}

\begin{remark}
 Since $B$ is bounded and nonnegative and $A^{-1}$ is bounded, $B A^{-1} B = 0$ is in fact equivalent to $A^{-1} B$ being
 nilpotent. Indeed, if we have $(A^{-1}B)^n = 0$ for some $n \in \NN$, $n \ge 2$, then
 $\spec( ( B^{1/2}A^{-1}B^{1/2} )^{n}) = \{ 0 \}$ since
 $\spec( ( B^{1/2}A^{-1}B^{1/2} )^{n}) \setminus \{ 0 \} = \spec ((A^{-1}B)^{n}) \setminus \{ 0 \} = \emptyset$.
 Furthermore, $BA^{-1}B = 0$ holds if and only if the map $z \mapsto (A^{-1}-zA^{-1}BA^{-1})B$ defines a~\emph{pseudo-resolvent};
 see~\cite{HilleP-57} for an introduction to this notion.
 In this case, this map also agrees with the pseudo-resolvent $(A + z B)^{-1} B$ belonging to the operator pencil.
\end{remark}

\begin{corollary}\label{cor1}
 If the equivalence in Proposition~\ref{prop} takes place with $B \neq 0$, then $A^{-1}BA^{-1} \neq 0$ and
 \begin{equation*}
  \frac{1}{\abs{t}\norm{A^{-1}BA^{-1}}+\norm{A^{-1}}}
  \le
  \dist\bigl(0,\spec(A+tB)\bigr)
  \le
  \frac{1}{\abs{t}\norm{A^{-1}BA^{-1}}-\norm{A^{-1}}}
 \end{equation*}
 for $t \in \RR$ with $\abs{t} > \norm{A^{-1}} / \norm{A^{-1}BA^{-1}}$.

 \begin{proof}
  First note that $A^{-1}BA^{-1} \neq 0$ is a consequence of the boundedness of $B$ and the fact that $A^{-1}$ is bijective as a
  map from the whole Hilbert space to the dense subspace $\Dom(A)$. Now, from identity~\eqref{eq:inverse} in Proposition~\ref{prop}
  we conclude that
  \begin{equation*}
   \abs{t}\norm{A^{-1}BA^{-1}}-\norm{A^{-1}} \le \norm{(A+tB)^{-1}} \le \abs{t}\norm{A^{-1}BA^{-1}} + \norm{A^{-1}}.
  \end{equation*}
  Since $\norm{(A+tB)^{-1}}^{-1}=\dist(0,\spec(A+tB))$, this proves the claim.
 \end{proof}%
\end{corollary}

Statement~(ii) in Proposition~\ref{prop} indicates how to construct examples where the set $\bigcup_{t \in \RR}\spec(A+tB)$ is not
the whole real line. The simplest is the following pedestrian's example of $2 \times 2$ matrices; there will be more sophisticated
examples below.

\begin{example}\label{example_2x2}
 Choosing
 \begin{equation*}
  A = \begin{pmatrix} 1 & 0\\ 0 & -1 \end{pmatrix},
  \quad
  B = \frac{1}{2}\begin{pmatrix} 1 & 1\\ 1 & 1 \end{pmatrix},
 \end{equation*}
 we see that $BA^{-1}B = 0$, whence $\bigcup_{t \in \RR}\spec(A+tB) \subset \RR \setminus \{ 0 \}$ by Proposition~\ref{prop}. In
 fact,
 \begin{equation*}
  \spec(A+tB) = \Bigl\{ -\frac{t}{2} \pm \sqrt{\frac{t^2}{4}+1} \Bigr\}\,,\quad t \in \RR,
 \end{equation*}
 so that $\bigcup_{t \in \RR}\spec(A+tB) = \RR \setminus \{ 0 \}$. Note that the distance of $0$ to the spectrum of $A+tB$ is given
 by $\dist(0,\spec(A+tB)) = (\abs{t}/2 + \sqrt{t^2/4+1})^{-1}$ and, therefore, behaves asymptotically for $\abs{t}\to\infty$
 exactly as predicted by Corollary~\ref{cor1}.
\end{example}

\begin{remark}
 The matrices of Example~\ref{example_2x2} appear in the context of critically damped linear systems, see for instance Example~9.3
 in~\cite{Ves} with the nilpotent matrix 
 \begin{equation*}
  A^{-1} B
  =
  \begin{pmatrix} 1 & 1\\ -1 & -1 \end{pmatrix}.
 \end{equation*}
\end{remark}

Clearly, the spectral point $0$ plays no particular role in the above considerations as we may replace $A$ by $A - \lambda$ for any
$\lambda \in \RR$. 
Since by Corollary~\ref{cor1} the distance of every protected point to the spectrum tends to zero in the
$\abs{ t } \to \infty$ limit, we deduce the following.

\begin{corollary}
 \label{cor2}
 Let $A$ and $B$ be as in Proposition~\ref{prop}.
 If $B \neq 0$, then the union $\bigcup_{t\in\RR} \spec(A+tB)$ is a dense subset of $\RR$.
\end{corollary}

The fact that the set $\bigcup_{t\in\RR} \spec(A+tB)$ is dense can, in fact, be refined to the following stronger statement.

\begin{lemma}\label{lem:count}
 Let $A$ and $B$ be as in Corollary~\ref{cor} with $B \neq 0$. Then, the set $\RR\setminus\bigcup_{t\in\RR}\spec(A+tB)$ is at most
 countable.

 \begin{proof}
  Pick any $y := Bx \neq 0$, and consider the function $f \colon \rho(A) \to \RR$ defined by
  \begin{equation}\label{eq:deff}
   f(z):=\langle x,B(A-z)^{-1}Bx \rangle = \langle y,(A-z)^{-1}y \rangle.
  \end{equation}
  For every $\lambda \in \RR$ with $B(A-\lambda)^{-1}B=0$ we find $f(\lambda)=0$. Hence, by Proposition~\ref{prop}, every point in
  $\RR \setminus \bigcup_{t\in\RR} \spec(A+tB)$ must be a (real) root of $f$. But $f$ can have at most countably many real roots.
  Indeed, we have $f'(z)=\langle y,(A-z)^{-2}y \rangle$ for all $z\in\rho(A)$, whence, in particular,
  \begin{equation*}
   f'(\lambda) = \norm{(A-\lambda)^{-1}y}^2>0 \quad\text{ for }\quad \lambda\in\rho(A) \cap \RR\,.
  \end{equation*}
  Thus, $f$ is strictly monotone on every interval of $\rho(A) \cap \RR$ and can therefore have at most one root there. Since
  $\rho(A) \cap \RR$ is an open subset of $\RR$, that is, an at most countable union of disjoint open intervals, this proves the
  claim.
 \end{proof}%
\end{lemma}

% \pagebreak
\begin{lemma}
 \label{lem:accumulation_points}
 Let $A$ and $B$ be as in Corollary~\ref{cor} with $B \neq 0$. If $\lambda$ is an accumulation point of
 $\RR \setminus \bigcup_{t\in\RR} \spec(A+tB)$, then $\lambda \in \spec(A+tB)$ for all $t \in \RR$. In particular,
 $\RR \setminus \bigcup_{t\in\RR} \spec(A+tB)$ contains none of its accumulation points, if any.

 \begin{proof}
  Let $(\lambda_k)$ be a sequence in $\RR \setminus \bigcup_{t\in\RR} \spec(A+tB)$ with $\lambda_k \to \lambda \in \RR$ as
  $k \to \infty$ and $\lambda_k \neq \lambda$ for all $k$, and assume that $\lambda \in \rho(A + t_0B)$ for some $t_0 \in \RR$.
  Set $\tilde{A} := A + t_0B$. Then, $(\lambda_k)$ is a sequence in $\RR \setminus \bigcup_{t\in\RR} \spec(\tilde{A}+tB)$, so that
  $B(\tilde{A}-\lambda_k)^{-1}B = 0$ for all $k$ by Proposition~\ref{prop}. Thus,
  the mapping $z \mapsto B(\tilde{A}-z)^{-1}B$ is analytic in $\lambda$ and has zeros at every $\lambda_k$, hence vanishes in a
  (complex) neighbourhood of $\lambda$. Therefore, again by Proposition~\ref{prop}, a real neighbourhood of $\lambda$ belongs to
  $\RR \setminus \bigcup_{t\in\RR} \spec(\tilde{A}+tB)$, so that the union
  $\bigcup_{t\in\RR} \spec(\tilde{A}+tB) = \bigcup_{t\in\RR} \spec(A+tB)$ is not dense in $\RR$. The latter is a contradiction to
  Corollary~\ref{cor2}.
 \end{proof}%
\end{lemma}

We now turn to the second assertion in Theorem~\ref{thm}. Here, let us first observe that the sets
$\RR\setminus\bigcup_{t\in\RR}\spec(A+tB)$ can indeed consist of every finite or countably infinite number of points. This can
easily be observed by inverting the construction in the proof of Lemma~\ref{lem:count} above:

\begin{example}
 \label{example_poles}
 Let $I$ be a finite or countably infinite set. On $\ell^2(I)$ consider the diagonal operator $A = \diag(\{\mu_k\}_{k \in I})$ with
 distinct real numbers $\mu_k$, and let $y \in \ell^2(I)$ be normalized with every entry non-zero. Finally, consider the rank-one
 projection $B:=\langle y,\cdot\rangle y\neq 0$. Clearly, $B(A-\lambda)^{-1}B = 0$ if and only if $\lambda$ is a zero of the
 function $f$ defined as in~\eqref{eq:deff} with $x := y$. Hence, by means of Proposition~\ref{prop}, the set
 $\RR \setminus \bigcup_{t\in\RR} \spec(A+tB)$ consists of all real zeros of $f$. Now, the mapping
 $\rho(A) \cap \RR \ni \lambda \mapsto f(\lambda)$ has poles in every $\mu_k$ and is strictly monotone on the intervals between
 them, cf.~the proof of Lemma~\ref{lem:count}. Therefore, it must have exactly one root on every such interval. We have thus found
 operators $A$ and $B$ with any finite or countably infinite number of ``protected'' points.
\end{example}

If we choose the $\mu_k$ in Example~\ref{example_poles} to have an accumulation point, then also the protected points between them
must have an accumulation point in $\RR$. However, by Lemma~\ref{lem:accumulation_points}, this accumulation point will not belong
to the set of protected points, but rather to each spectrum $\spec(A+tB)$, $t \in \RR$.
 
The method of Example~\ref{example_poles} does not allow us to directly choose the protected points.
This can be achieved by the following construction, which concludes the proof of Theorem~\ref{thm}:

Let $P \subset \RR$ be a finite or countably infinite discrete set. Choose an orthonormal basis
$(\psi_\lambda)_{\lambda\in P}$ of $\cH = \ell^2(P)$, and consider the (not necessarily bounded) self-adjoint operator $K$ on $\cH$
with $K \psi_\lambda = \lambda \psi_\lambda$ for all $\lambda \in P$.
Clearly, every $\lambda \in P$ is an isolated eigenvalue of multiplicity one. Thus, we have $\spec(K) = \overline{P}$, and the
operator $K$ has simple spectrum (see, e.g.,~\cite[Section~5.4]{Schm12} for a definition). In particular,
by~\cite[Proposition 5.20]{Schm12} there exists a~\emph{cyclic vector} $v \in \cH$, which means that $\{ v, K v, K^2v, \dots \}$
spans a dense subspace of $\cH$.

We now add one extra dimension and define self-adjoint operators $A$ and $B$ on the Hilbert space $\tilde \cH := \cH \oplus \CC$ by 
\begin{equation}
 \label{eq:defABwithK}
 A
 :=
 \begin{pmatrix}
  K & v\\
  \langle v,\cdot \rangle & 0 \\
 \end{pmatrix}
 ,
 \quad
 B
 :=
 \begin{pmatrix}
  0 & 0\\
  0 & 1\\
 \end{pmatrix}
 .
\end{equation}

\begin{lemma}
 \label{lem:cyclic}
 For $A$ and $B$ as in~\eqref{eq:defABwithK} we have $\bigcup_{t \in \RR} \spec(A+tB) = \RR \setminus P$.
 
 \begin{proof}
  For each $t \in \RR$ the operator $A + t B$ is a finite rank perturbation of $K \oplus 0$, so that its spectrum consists of the
  essential spectrum of $K$, that is, the accumulation points of $P$, and isolated eigenvalues of finite multiplicity.
  
  Let $\lambda \in P$. We need to show that $\lambda \in \rho(A + t B)$ for all $t \in \RR$. Here, it suffices to see that each
  $A + t B - \lambda$ has trivial kernel since $\lambda$ is not in the essential spectrum of $A + t B$. To this end, let
  $x \oplus \alpha \in \Ker(A + tB - \lambda)$, that is,
  \begin{equation}
   \label{eq:kernel}
   \begin{pmatrix}
    (K - \lambda)x + \alpha v\\
    \langle v,x \rangle + (t - \lambda)\alpha
   \end{pmatrix}
   =
   \begin{pmatrix}
    0\\
    0
   \end{pmatrix}
   .
  \end{equation}
  We expand $v = \sum_\lambda \beta_\lambda \psi_\lambda$ in the basis $(\psi_\lambda)_\lambda$ and note that cyclicity of $v$
  forces all $\beta_\lambda$ to be non-zero. Now, by~\eqref{eq:kernel} we have $(K - \lambda) x = -\alpha v$. Since
  $(K - \lambda) x$ is orthogonal to $\psi_\lambda$ and $\beta_\lambda \neq 0$, this implies that $\alpha = 0$. Hence,
  $(K-\lambda)x = 0$ and, therefore, $x = \langle \psi_\lambda,x \rangle \psi_\lambda$ is a multiple of $\psi_\lambda$. It then
  follows from~\eqref{eq:kernel} that $0 = \langle v,x \rangle = \overline{\beta_\lambda} \langle \psi_\lambda,x \rangle$, which
  yields that also $x = 0$.

  Conversely, let $\lambda \in \RR \setminus \bigcup_{t\in\RR} \spec(A + tB)$. Since, in particular, $\lambda \in \rho(A)$, there
  is $x \oplus \alpha \in \tilde{\cH} \setminus \{0\}$ with $(A - \lambda)(x \oplus \alpha) = 0 \oplus 1$. Now, by
  Proposition~\ref{prop} we have $B(A - \lambda)^{-1}B = 0$ and, therefore,
  $0 = B(A - \lambda)^{-1}B(0 \oplus 1) = 0 \oplus \alpha$. This implies that $\alpha = 0$ and, in turn, $(K - \lambda)x = 0$ with
  $x \neq 0$. Hence, $\lambda$ is an eigenvalue of $K$, that is, $\lambda \in P$.
 \end{proof}%
\end{lemma}

\subsubsection*{Acknowledgments}
The authors are grateful to the anonymous referees for helpful remarks on the manuscript. M.T.~was supported in part by the
European Research Council starting grant 639305 (SPECTRUM) and is indebted to Sabine B\"ogli and Sasha Sodin for valuable
discussions.

\end{document}